\documentclass[11pt,a4,english]{amsart}

\usepackage{amsmath}   
\usepackage{amsfonts,amssymb}

\usepackage{graphicx} 
\usepackage{url}      

\usepackage{amsfonts}
\usepackage{babel}
\usepackage{latexsym}
\usepackage{verbatim}
\usepackage{color}

\newcommand\dela[1]{}

\numberwithin{equation}{section}

\newtheorem{theorem}{Theorem}[section]
\newtheorem{assumption}[theorem]{Assumption}

\newtheorem{lemma}[theorem]{Lemma}
\newtheorem{corollary}[theorem]{Corollary}

\newtheorem{definition}[theorem]{Definition}
\newtheorem{remark}[theorem]{Remark}
\newtheorem{proposition}[theorem]{Proposition}

\newcommand{\fh}{{\mathfrak{h}}}
\newcommand{\fs}{{\mathfrak{s}}}

\newcommand{\fhs}{\mathfrak{h} ^ s}

\newcommand{\bcase}{\begin{cases}}
\newcommand{\ecase}{\end{cases}}

\newcommand\del[1]{}

\newcommand\cadlag{c{\`a}dl{\`a}g }

\newcounter{gg11}
\newenvironment{list-a}
{\begin{list} {{\rm (\alph{gg11})}}
{\usecounter{gg11}
\setlength{\leftmargin}{0.3cm}
\setlength{\topsep}{0.1cm}
\setlength{\itemsep}{0.0cm}
\setlength{\parsep}{0.1cm}
\setlength{\itemindent}{0.0cm}
\setlength{\parskip}{0.0cm}}}
{\end{list}}
\newcounter{gg111}
\newenvironment{list-n}
{\begin{list} {{\rm (\roman{gg111})}}
{\usecounter{gg111}
\setlength{\leftmargin}{0.5cm}
\setlength{\topsep}{0.1cm}
\setlength{\itemsep}{0.0cm}
\setlength{\parsep}{0.1cm}
\setlength{\itemindent}{0.0cm}
\setlength{\parskip}{0.0cm}}}
{\end{list}}
\newcounter{gg1}

\newcounter{lil}

\newcommand{\CMM}{{{M}_I([0,T] \times Z)}}
\newcommand{\CMME}{{{M}_I([0,1]\times Z )}}

\newcommand{\CMR}{{M_+(Z)}}
\newcommand{\lk}{\left}
\newcommand{\lqq}{\lefteqn}
\newcommand{\rk}{\right}

\newcommand{\LL}{{\rm I \kern -0.2em L}}

\newcommand{\ep} {\varepsilon }

\newcommand{\be} {\begin{enumerate} }
\newcommand{\ee} {\end{enumerate} }

\newcommand{\CZ}{{{ \mathcal Z }}}

\newcommand{\CX}{{{ \mathcal X }}}

\newcommand{\CLaw}{{{ \mathcal Law }}}

\newcommand{\CB}{{{ \mathcal B }}}
\newcommand{\CM}{{{ \mathcal M }}}

\newcommand{\CF}{{{ \mathcal F }}}

\newcommand{\CL}{{{ \mathcal L }}}

\newcommand{\RR}{{\mathbb{R}}}

\newcommand{\DD}{{\rm I \kern -0.2em D}}
\newcommand{\dd}{{\rm I \kern -0.2em D}}
\newcommand{\bNN}{{\bar{ \mathbb{N}}}}

\newcommand{\PP}{{\mathbb{P}}}

\newcommand{\TT}{{\rm I \kern -0.2em T}}

\newcommand{\DEQS}{\begin{eqnarray*}}
\newcommand{\EEQS}{\end{eqnarray*}}
\newcommand{\DEQSZ}{\begin{eqnarray}}
\newcommand{\EEQSZ}{\end{eqnarray}}
\newcommand{\DEQ}{\begin{eqnarray}}
\newcommand{\EEQ}{\end{eqnarray}}

\newcommand{\ttin}{0\le t<\infty}

\newcommand{\e}{\mathbb{E}}
\newcommand{\p}{\mathbb{P}}

\newcommand{\fmath}{\mathbb{F}}

\newcommand{\bcal}{\mathcal{B}}
\newcommand{\fcal}{\mathcal{F}}
\newcommand{\zcal}{\mathcal{Z}}

\begin{document}

\title[Uniqueness in Law of the stochastic convolution \today]{
Uniqueness in Law of the stochastic convolution process driven by L\'evy noise}
\thanks{This work was supported by the FWF-Project P17273-N12}
\author[Z. Brze{\'z}niak,
E. Hausenblas and E. Motyl ]{Zdzis{\l}aw Brze{\'z}niak,
Erika Hausenblas and El\.zbieta Motyl}
\address{Department of Mathematics \\
 University of York, Heslington, York YO10
5DD, UK} \email{zdzislaw.brzezniak@york.ac.uk}
\address{Department of Mathematics and Information Technology, Montanuniversity of Leoben,
Franz Josef Strasse 18, 8700 Leoben, Austria}
\email{erika.hausenblas@leoben.ac.at}
\address{Department of Mathematics and Computer Sciences, University of \L\'{o}d\'{z}, ul. Banacha 22,
91-238 \L \'{o}d\'{z}, Poland}
\email{emotyl@math.uni.lodz.pl}

\date{\today}

\begin{abstract}
We will give a proof of the following fact.
If $\mathfrak{A}_1$ and $\mathfrak{A}_2$, $\tilde \eta_1$ and $\tilde \eta_2$, $\xi_1$ and $\xi_2$ are two examples of filtered probability spaces, time homogeneous compensated Poisson random measures,  and progressively measurable Banach space valued processes such that the laws on $L^p([0,T],{L}^{p}(Z,\nu ;E))\times \CM_I([0,T]\times Z)$ of the pairs  $(\xi_1,\eta_1)$ and $(\xi_2,\eta_2)$ 
  are equal, and $u_1$ and $u_2$ are the corresponding stochastic convolution processes, then the laws on $ (\DD([0,T];X)\cap L^p([0,T];B)) \times L^p([0,T],{L}^{p}(Z,\nu ;E))\times \CM_I([0,T]\times Z) $, where $B \subset E \subset X$, of the triples  $(u_i,\xi_i,\eta_i)$, $i=1,2$, are equal as well. By $\DD([0,T];X)$ we denote the
  Skorokhod space of $X$-valued processes.

\end{abstract}

\maketitle
\baselineskip 16pt \markright{Uniqueness in Law of the stochastic convolution}

\renewcommand{\labelenumi}{\alph{enumi}.)}

\noindent
\textbf{Keywords and phrases:} {Poisson random measure, stochastic convolution process, uniqueness in law,
stochastic partial differential equations.}

\bigskip \noindent
\textbf{AMS subject classification (2002):} {Primary 60H15;
Secondary 60G57.}

\section{Introduction}\label{sec-intro}

A solution to a stochastic partial differential equation (SPDE) driven by a  martingale
can be represented as a stochastic convolution process with respect to that martingale.
If the coefficients of this SPDE are globally Lipschitz continuous,  the solution can be found
by an application of the Banach Fixed Point Theorem, see for instance  Da Prato and Zabczyk \cite{DaPrato+Zabczyk_1992}, Brze\'{z}niak \cite{Brz_1997},  Hausenblas \cite{levy2} or St. Loubert Bi\'e \cite{980.39765}.
However, if the coefficients  are  only continuous, a typical approach to proof the existence of
a solution is via the  use the Skorokhod Theorem on representation of a weakly convergent sequences of measures by a.s. convergent random variables.  But then
the  original probability space becomes lost and the solution is  defined on a new probability space.
For this reason it is important to know whether
the joint law of the triplet consisting of the driving martingale, the integrand
 and the stochastic convolution process, remains the same.
If the driving martingale is  continuous  one can use
various versions of the martingale representation  theorem to justify
this equality, see for instance Da Prato and Zabczyk \cite[Section 8]{DaPrato+Zabczyk_1992} and the references therein.
If the driving martingale  is  purely discontinuous we have been unable to find an appropriate
 embedding theorem to proof the required equality. The aim of the present paper
is to fill this gap in the existing literature. To be precise,
we will show that changing the underlying probability space without changing the laws of the integrand and the driving martingale does not lead to a change of the law of the corresponding  triplet (consisting of the driving martingale, the integrand and the stochastic convolution process).

We encountered this  sort of a difficulty  in our recent paper \cite{reactdiff} in which we studied the existence of a solution to a stochastic reaction diffusion equation, with only continuous coefficient,  driven by a purely  discontinuous L\'evy process. It turned out that the use of the result presented in this  paper was essential. We believe that this result is  interesting on its own (as well as it will be  applicable in other situations) and, hence,  we have decided to publish it separately.

Let us stress here that the difficulty lies in the problem of the existence of a c\`{a}dl\`{a}g modification of the stochastic convolution processes. Let us notice that this is not always true as has been recently shown by Brze{\'z}niak et al \cite{Brzezniak_noncadlag_2010}.
However, in  paper  \cite{Brz+Haus+Zhu_2012} the authors give certain positive answer to this problem in the case when a martingale type Banach space satisfies some additional rather restrictive assumption.
The crucial point is the proof of appropriate  maximal inequalities for stochastic convolution.
In the general case, this problem seems to be open.

We finish this introduction with a brief summary of our results. But for that aim we need to introduce the basic notation we will be using throughout the whole paper.

\bigskip  \noindent
\textbf{Notation:}
By $\mathbb{N}$ we denote the set of natural numbers (including $0$) and by $\bar{\mathbb{N}}$ we denote the set $\mathbb{N}\cup\{\infty\}$.
For a measurable space $(Z,\CZ)$ by $M_I(Z) $ we denote the family of all $\bar{\mathbb{N}}$-valued measures on $(Z,\CZ)$ and by $\CM_I(Z) $ the $\sigma$-field on $M_I(Z)$ generated by functions
$i_B:{M}_{I}(Z)\ni\mu \mapsto \mu(B)\in \bNN$, $B\in \CZ$.
By $\CMR$ we denote the set of all non negative and $\sigma$-finite  measures on $(Z,\CZ )$.

By $-A$ we denote an infinitesimal generator of a $C_0$ semigroup  $\big(S(t)\big)_{\ttin   }$  on a Banach space $(E, {\vert \cdot \vert }_{E} )$. We assume that $p\in (1,2]$ is fixed and that  $E$ is a separable martingale type $p$ Banach space, see  for instance   \cite[Appendix A]{Brz+Haus_2009}.

We assume that  $(Z,\CZ)$ is a measurable space,  $\mathfrak{A}=(\Omega,\CF,\mathbb{F},\PP)$, where $\mathbb{F}=(\CF_t)_{\ttin   }$, is a filtered probability space and
$\tilde{\eta }$ is a  time homogeneous compensated Poisson random measure on $(Z,\CZ)$ over  $\mathfrak{A}$\del{,
with intensity measure $\nu$}.
Let $\xi:  \RR_+ \times \Omega \times Z\to E$ be a progressively measurable  process such that
 for every $T >0 $
\begin{equation*} \label{cond-1.01}
\mathbb{E} \int_0^T \int_Z\vert \xi(r,z)\vert_E^p\,
\nu(dz)\,dr < \infty,
\end{equation*}
where $\nu$ is the intensity measure of $\tilde{\eta }$.
We consider the following It\^o type SPDE in the space $E$
\begin{equation} \label{eqn-Langevin-Poisson}
\quad\quad\left\{\begin{array}{ll}
  & d u(t) + Au(t) \; dt
 = \int_{Z}\xi(t;z)\, \tilde \eta(dt;dz),\quad t > 0,
\\ & u(0) = 0.
\end{array}
\right.
\end{equation}
A solution to problem  \eqref{eqn-Langevin-Poisson} can be defined  to be the following stochastic convolution process with respect  to $\tilde \eta$:
\begin{equation} \label{eqn-SC-Poisson}
u(t) := \int_0 ^ t \int_{Z}  S(t-r)\, \xi(r,z)\: \tilde{\eta }(dr;dz), \quad  t>0.
\end{equation}

In general, we can expect that $u$ has a \cadlag modification provided that it is considered as a process in an appropriately chosen Banach space $X$ such that $E \subset X$. Moreover, by smoothing properties of some types of semigroups $u$, considered as a process in some smaller Banach space $B \subset E$, has
$p$-integrable paths. So, we formulate our results in the general framework. However, we present also applications and give examples of the spaces $X$ and $B$.

The aim of this paper is to prove the following result.
If $\mathfrak{A}_1$ and $\mathfrak{A}_2$, $\tilde \eta_1$ and $\tilde \eta_2$, $\xi_1$ and $\xi_2$ are two examples of filtered probability spaces, time homogeneous compensated Poisson random measures,  and progressively measurable processes such that the laws on
 $L^p([0,T],{L}^{p}(Z,\nu ;E))\times {M}_{I}([0,T]\times Z)$ of the pairs  $(\xi_1,\eta_1)$ and $(\xi_2,\eta_2)$
  are equal, and $u_1$ and $u_2$ are the corresponding stochastic convolution processes, then the laws on $(\DD([0,T];X)\cap L^p(0,T;B) \times  L^p([0,T],{L}^{p}(Z,\nu ;E))\times {M}_{I}([0,T]\times Z))$ of the triplets  $({u}_{i},\xi_i,\eta_i)$, $i=1,2$, are equal as well. Here, by $\DD([0,T];X)$ we denote the
  Skorokhod space of $X$-valued processes.
The case of stochastic integrals, i.e. the stochastic convolution processes with $A=0$,  was studied in  \cite{Brz+Haus_2009-unique}.

\bigskip
This paper is organised as follows. Section \ref{sec-prel} contains some probabilistic preliminaries.
The main results are stated in Section \ref{sec-main}.
In Section \ref{sec-appl} we present applications of the general results to stochastic differential equations.
 Auxiliary results on the Haar projection in  the space ${L}^{p}([0,T];E)$ and the dyadic projection in the Skorokhod space are contained in Appendices B and C, respectively.

\subsection*{Acknowledgments} The authours would like to thank an anonymous referee for reading the paper twice. His comment, questions  and suggestions have lead to a substantial improvement of manner of  presentation of the paper.

\section{Preliminaries} \label{sec-prel}

Let us first recall a definition of time homogeneous Poisson random measures over a filtered probability space.
\begin{definition}\label{def-Prm}
Let $(Z,\zcal )$ be a measurable space, $\nu \in {M}_{+}(Z)$ and let $\mathfrak{A}=(\Omega,\fcal ,\fmath ,\PP)$  be a
filtered probability space with the filtration $\fmath = {({\fcal }_{t})}_{0\le t < \infty }$.
A {\sl time homogeneous Poisson random measure \it $\eta $  on $(Z,\zcal )$  with intensity measure $\nu $}
 \it over $\mathfrak{A}$ is a measurable map $$\eta : (\Omega , \fcal ) \to ({M}_{I}({\mathbb{R}}_{+}\times Z), {\mathcal{M}}_{I}({\mathbb{R}}_{+}\times Z))$$
satisfying the following conditions
\begin{trivlist}
\item[(i)] for all $B \in \bcal ({\mathbb{R}}_{+}) \otimes \zcal $, $\eta (B):= {i}_{B} \circ \eta : \Omega \to \overline{\mathbb{N}} $ is a Poisson random measure with parameter
 $\e [\eta (B)]$,
\item[(ii)] $\eta $ is independently scattered, i.e. if the sets ${B}_{j} \in \bcal ({\mathbb{R}}_{+}) \otimes \zcal $, $j=1,...,n$, are disjoint then the random variables $\eta ({B}_{j})$, $j=1,...,n$, are independent,
\item[(iii)] for all $U \in \zcal $ and $I \in \mathcal{B}({\mathbb{R}}_{+})$
$$
    \e [\eta (I \times U)] = \lambda (I) \nu (U) ,
$$
where $\lambda $  is the Lebesgue measure,
\item[(iv)] for all $U \in \zcal $  the
$\overline{\mathbb{N}}$-valued process ${(N(t,U))}_{t \ge 0 }$ defined by
$$
   N(t,U):= \eta ((0,t]\times U) , \qquad t \ge 0 ,
$$
is $\fmath $-adapted and its increments are independent of the past, i.e. if $t>s \ge 0 $, then
$N(t,U) - N(s,U)= \eta ((s,t]\times U)$  is independent of ${\fcal }_{s}$.
\end{trivlist}
\end{definition}

Let $\eta $ be a time homogeneous Poisson random measure with intensity $\nu \in {M}_{+}(Z)$ over
$\mathfrak{A}$.
We will denote by $\tilde{\eta }$ the \it  compensated Poisson random measure \rm  defined by
$\tilde{\eta }: = \eta - \gamma $, where the compensator
 $\gamma : \bcal ({\mathbb{R}}_{+})\times \zcal  \to {\mathbb{R}}_{+}$ satisfies in our case the following equality
$$
     \gamma (I \times A) = \lambda (I) \nu (A) ,
\qquad I \in \bcal ({\mathbb{R}}_{+}) , \quad  A \in \zcal .
$$

\bigskip

Let us now state an assumption we will be using throughout the whole paper.
\begin{assumption}\label{assumption-main}
We suppose that $\big(E,\vert \cdot\vert_E\big)$ is a separable Banach space of martingale type $p$,  where   $1< p\le 2$.
 \end{assumption}

 In  \cite{Brz+Haus_2009} there is proven that there exists a unique
continuous linear operator $I$ which associates   to  each progressively measurable
process $\xi: {\mathbb{R}}_{+} \times \Omega \times Z \to E$  with
\begin{equation} \label{cond-2.01}
 \e \Bigl[ \int_{0}^{T} \int_{Z} \vert \xi (r,x) \vert_E^{p} \, \nu (dx) dr  \Bigr] < \infty
\end{equation}
for every $T>0$, an adapted $E$-valued c\`{a}dl\`{a}g process
$$
  {I}_{\xi , \tilde{\eta }} (t):=\int_{0}^{t} \int_{Z} \xi (r,x) \tilde{\eta } (dr,dx), \quad t \ge 0
$$
such that if a process $\xi$ satisfying  the above condition (\ref{cond-2.01}) is a
random step process with representation
$$
\xi(r,x) = \sum_{j=1} ^n 1_{(t_{j-1}, t_{j}]}(r) \, \xi_j(x),\quad x\in Z,\quad  r\ge 0,
$$
where $\{t_0=0<t_1<\ldots<t_n<\infty\}$ is a finite partition of $[0,\infty)$ and
for all $j\in \{ 1, \ldots ,n \} $,    $\xi_j$ is  an $E$-valued $\CF_{t_{j-1}}$--measurable $p$-summable simple  random variable,
then
\begin{equation} \label{eqn-2.02}
 I_{\xi,\tilde{\eta }} (t) 
=\sum_{j=1}^n  \int_Z  \xi_j (x) \,\tilde \eta \lk((t_{j-1}\wedge t, t_{j}\wedge t] ,dx \rk),\quad   t\ge 0.
\end{equation}

In the recent paper  \cite{Brz+Haus_2009-unique} there is shown that this
continuous linear operator is unique in a weak sense. In particular,
there is proved  that for every $T>0$
the law of the triplet $(I_{\xi,\tilde{\eta }},\xi,\eta)$ on
$$
 L^p([0,T],E)\cap\DD([0,T] ,E)\times L^p([0,T] ,L^p(Z,\nu;E))\times \CMM
$$
depends only on the law of the pair $(\xi,\eta)$ on the corresponding space.
\del{is independent from its representation on the underlying probability space. }
An  important point to mention here is that $\eta$ is a time homogeneous Poisson random measure and
that the filtration generated by $\xi$ is nonanticipative with respect to $\eta$.
In this paper, we will show that a similar results hold also for the stochastic convolution process
defined in  \eqref{eqn-SC-Poisson}.

\section{Main results}\label{sec-main}

We assume that  $(Z,\CZ)$ is a measurable space,  $\mathfrak{A}=(\Omega,\CF,\mathbb{F},\PP)$, where $\mathbb{F}=(\CF_t)_{\ttin   }$, is a filtered probability space and
$\tilde{\eta }$ is a  time homogeneous compensated Poisson random measure on $(Z,\CZ)$ over  $\mathfrak{A}$\del{,
with intensity measure $\nu$}.
Let us fix $1<p\le 2$.
Suppose that $(E,\vert \cdot \vert_E )$ is a separable Banach space of martingale type $p$.

\subsection{Stochasic convolution with a ${\mathcal{C} }_{0}$-semigroup}  \label{subsec-C_0-semigroup}

The following will be a standing assumption in most of this section.
\begin{assumption}\label{assumption-semigroup}
The operator $-A$ is an infinitesimal generator of  a ${\mathcal{C} }_{0}$-semigroup  ${(S(t))}_{0\le t < \infty }$  on the space $E$.
\end{assumption}

Let $\xi: \Omega\times \RR_+\times Z\to E$ be a progressively measurable  process such that for all $T>0$
\begin{equation} \label{E:Main_cond-1.01}
\mathbb{E} \int_0^T \int_Z\vert \xi(r,z)\vert_E^p\,
\nu(dz)\,dr < \infty,
\end{equation}
where $\nu$ is the intensity measure of $\tilde{\eta }$.
Let us consider the stochastic convolution process
$u=u_{\xi,\tilde{\eta}} $ defined by
\begin{equation}  \label{E:Main_stoch_conv}
 u(t) = \int_0 ^ t \int_Z S(t-s)
\xi(s,z)\, \tilde \eta(ds,dz), \quad t\ge 0.
\end{equation}

\begin{remark} \label{R:Stoch_conv_L^p}
Note that $\{ u(t), t \ge 0 \} $ given in \eqref{E:Main_stoch_conv} is well defined process such that for every $T>0$,  $\p $-a.s.,
$u \in {L}^{p}([0,T];E)$.
\end{remark}

\bigskip
Let us  recall the following elementary definition.

\begin{definition}\label{def-equality}  
Let $(X,\CX)$ be a measurable space.
When we say that $\xi_1$ and $\xi_2$ have the same law on $X$ (and write $\CL aw(\xi_1)=\CL aw(\xi_2)$ on $X$), we mean that $\xi_i$, $i=1,2$, are $X$-valued random variables defined over some
probability spaces $(\Omega_i,\CF_i,\PP_i)$, $i=1,2$, such that
$$
\PP_1 \bar{\circ}\, \xi_1= \PP_2\bar{\circ}\, \xi_2,
$$
where $\PP_i\bar{\circ}\, \xi_i(A) = \PP_i(\xi ^ {-1}_i(A))$, $A\in\CX$, $i=1,2$, is a   probability measure on $(X,\CX)$ called the law of $\xi_i$.
\end{definition}

\bigskip  \noindent
\bf Remark. \rm (see \cite[Theorem II. 3.2]{para})
If $X$ is a separable metric space then $\PP_i\bar{\circ }\, \xi_i $, $i=1,2$, are Radon measures.

\bigskip
Now, we are ready to  state our main result. We will formulate it in the most general form possible. In the next section we will present a couple of important applications.

\begin{theorem} \label{thm-main}
Let us assume that $(E,\vert \cdot \vert_E)$ be a separable  Banach space of martingale type $p$, where $1<p\le 2$. Let us also assume that  Assumption \ref{assumption-semigroup} is satisfied.
Suppose that for  $i=1,2$
\begin{trivlist}
\item[(a)] ${\mathfrak{A}}_{i}=(\Omega_i,\CF_i,{\mathbb{F}}_{i},\PP_i)$, where ${\mathbb{F}}_{i}= (\mathcal{F}^i_t)_{\ttin}$,
is a  complete filtered probability space;
\item[(b)] $\eta_i$
 is a time homogeneous Poisson random measure on $(Z,\CZ)$
with intensity measure $\nu$ over ${\mathfrak{A}}_{i}$;  
\item[(c)]
$\xi_i$ is a  progressively measurable process over ${\mathfrak{A}}_{i}$
 satisfying condition \eqref{E:Main_cond-1.01}.
\end{trivlist} 
Let $T>0$. We put for $i=1,2$
\begin{equation} \label{def_u_i}
{u}_{i} (t):= u_{\xi_i,{\tilde{\eta}}_{i}} (t) = \int_{0}^{t} \int_{Z} S(t-r)\,{\xi }_{i}(r;z) \, \tilde  {\eta }_{i}(dr,dz), \quad  t\in [0,T].
\end{equation}
Assume finally that    $\CL aw( (\xi  _1,\eta_1))=\CL aw((\xi _2,\eta_2)) $ on
$L ^ p([0,T];L ^ p(Z,\nu;E)) \times \CMM  $.
\begin{trivlist}
\item[(i)]  If  $X$ is a separable Banach space such that $E \subset X $ continuously
and the processes $u_i$, $i=1,2$,
have $X$-valued c\`{a}dl\`{a}g modifications   (denoted by the same symbols), then
$\CL aw( (u _1,\xi  _1,\eta_1))=\CL aw((u _2,\xi _2,\eta_2)) $ on
 \[  \mathbb{D}([0,T];X) \times  L^p([0,T] ;L ^ p(Z,\nu;E))\times \CMM.\]
\item[(ii)]
If moreover, $B$ is separable Banach spaces such that $B \subset E$ continuously and the processes $u_i$, $i=1,2$,  have ${\p }_{i}$-almost all paths in ${L}^{p}(0,T;B)$, then $\CL aw( (u _1,\xi  _1,\eta_1))=\CL aw((u _2,\xi _2,\eta_2)) $ on
\[ (\DD([0,T] ,X)\cap L^p([0,T] ;B))\times  L^p([0,T];L^p(Z,\nu ;E))\times \CMM .\]
\end{trivlist}
\end{theorem}

\bigskip
\begin{remark}
The claim that $\CL aw( (u _1,\xi  _1,\eta_1))=\CL aw((u _2,\xi _2,\eta_2)) $ on
$$(\DD([0,T] ,X)\cap L^p([0,T] ;B))\times  L^p([0,T];L^p(Z,\nu ;E))\times \CMM $$
is essentially stronger than a  similar claim that
$\CL aw( u _1)=\CL aw( u _2)$ on $\DD([0,T] ;X)\cap L^p([0,T] ;B)$.
\end{remark}

Before we embark on  with the proof of Theorem \ref{thm-main} let us introduce
some useful notation. Given a  filtered probability space $\mathfrak{A}=(\Omega,\CF,\mathbb{F},\PP)$,
 and  a Banach  space  $Y$,
by  $\mathcal{N}([0,T] \times\Omega ;Y)$ we denote the space of (equivalence classes of)
progressively measurable functions $\xi :[0,T] \times\Omega \to Y$.

For $ q\in (1,\infty ) $ we  set
\begin{equation}
\mathcal{M}^q([0,T] \times\Omega  ;Y)=
\Bigl\{ \xi \in \mathcal{N}([0,T] \times\Omega  ;Y):\mathbb{E}\int_
0^T\vert\xi (t)\vert_Y^q\,dt<\infty \Bigr\} .
\label{def:Mq}
\end{equation}
Let $\mathcal{N}_{\rm step}([0,T] \times\Omega   ;Y)$ be the space of all $
\xi \in \mathcal{N}([0,T] \times\Omega    ;Y)$ for which  there exists a partition
$0=t_0<t_1<\cdots <t_n=T$ such that for $k\in\{1,\cdots,n\}$ and for $t\in (t_{k-1},t_{k}]$,
$\xi (t)=\xi (t_k)$ is $\mathcal{F}_{t_{k-1}}$-measurable. We put $\mathcal{M}_{\rm  step}^q
=\mathcal{M}^q\cap \mathcal{N}_{\rm step}$.
Note that $\mathcal{M}^q([0,T] \times \Omega ;Y)$ is a closed subspace of $L^q([0,T] \times\Omega ;Y)
\cong L^q(\Omega ;L^q([0,T];Y))$.

\begin{proof}[Proof of  Theorem \ref{thm-main}]
Note that  the second assertion  is a consequence  the first one.  Indeed,  by the first claim,  the triplets $(u_{\xi_i,{\tilde{\eta}}_i},\xi_i,\eta_i)$, $i=1,2$,  have the same law
on $\mathbb{D}([0,T];X)\times L^p([0,T];{L}^{p}(Z,\nu ;E))\times M_I([0,T]\times Z)$.
Since $\DD([0,T];X) \cap L^p(0,T;B) \hookrightarrow \mathbb{D}([0,T];X)$ continuously, in view of Proposition \ref{prop-main}
 the same  triplets
have equal laws
on the space
$$ \DD([0,T];X) \cap L^p(0,T;B) \times L^p([0,T];{L}^{p}(Z,\nu ;E))\times M_I([0,T]\times Z).$$
Thus it is enough to prove assertion (i).
In order to make the use of the {Haar projection  more transparent} we will assume in this
part of the paper that $T=1$.
Let us fix a  filtered probability space $\mathfrak{A}=(\Omega,\CF, \fmath ,\PP)$, where $\fmath = (\CF_t)_{ t\in [0,1]}$, and let us introduce the family of linear operators  $(\Phi_t)_{t \in [0,1]}$
from
$ L^p( [0,1]  ;L^p(Z,\nu;E))$ to $ L^p( [0,1]  ;L^p(Z,\nu;X))$, defined by the following formula
$$
(\Phi _t \xi )(s) := 1_{[0,t)}(s)  S(t-s) \xi (s) ,\quad
\xi \in   L^p([0,1] ;L^p(Z,\nu;E)) , \quad s\in [0,1]  .
$$
Note that for any $t\in [0,1]$, the operator $\Phi_t$  is well defined  bounded and linear.
Indeed,
there exists a constant $M>0$ such that for all $t \in [0,1]$, ${|S(t)|}_{\mathcal{L}(E)} \le M$.
Hence by the continuity of the embedding $\iota :E \hookrightarrow X$  we obtain for any $t\in [0,1]$
\begin{eqnarray}
\label{E:Phi_t_ineq}    \nonumber
& & \int_{0}^{1} \| ({\Phi }_{t}\xi)(s)\|_{{L}^{p}(Z,\nu ;X)}^{p}
   =  \int_{0}^{1} \| 1_{[0,t)}(s) \iota \circ S(t-s) \xi (s) \|_{{L}^{p}(Z,\nu ;X)}^{p} \, ds  \\
 &  & =  \int_{0}^{1} \int_{Z}\| 1_{[0,t)}(s) \iota \circ S(t-s) \xi (s)(z) \|_{X}^{p} \, d\nu (z) ds \\
 && \leq C{M}^{p}  \int_{0}^{1} \vert \xi(s)(z) \vert_{E}^{p} \, d\nu (z) ds =CM^p \Vert \xi\Vert^p_{L^p(0,1;L^p(Z,\nu ;E))}.
      \nonumber
\end{eqnarray}
Thus the operator ${\Phi }_{t}$ is well defined and bounded.

Notice that the stochastic convolution process ${u}_{\xi , \tilde{\eta }}$ defined in \eqref{E:Main_stoch_conv} can be expressed in terms of the map ${\Phi }_{t}$. Indeed, identifying a process $\xi $ with a map $\xi : \Omega \to {L}^{p}([0,1];{L}^{p}(Z,\nu ;E))$ we have the following equality
$$
u_{ \xi,\tilde \eta }(t) =  I_{\Phi_t \circ \xi ,\tilde \eta}(t)=
 \int_0 ^ t \int_Z (\Phi_t \circ \xi)(s,z)\,\tilde{\eta }(ds,dz),\quad t\in[0,1] .
$$
Going back to our original problem we assume  that  ${\mathfrak{A}}_{i}=(\Omega_i,\CF_i,{\fmath }_{i},\PP_i)$, $i=1,2$, where ${\fmath }_{i}= (\CF ^ i_t)_{\ttin} $,
are two fixed complete, filtered probability spaces and $\xi_i\in\CM^p([0,1] \times\Omega_i  ;L^p(Z,\nu,E))$, $i=1,2$.
For fixed $t\in [0,1]$ we approximate $\Phi_t \circ \xi_i$, $i=1,2$,
by a  sequence
$$\lk(\fhs _n \circ \Phi_t \circ \xi _i \rk)_{n\in\mathbb{N}}\subset \CM^p_{\rm step}([0,1] \times\Omega_i    ;L^p(Z,\nu,E)),\quad  i=1,2,
$$
where
$\fhs_n$ is the shifted Haar projection operator
in the space ${L}^{p}([0,1];{L}^{p}(Z,\nu ,X))$, defined in
(\ref{haar-projection}).
Let us note that by Proposition \ref{cont-haar-shift}-(i) the shifted Haar projection $\fh ^s _n$ is also a  continuous operator from $L^p([0,1];L^p(Z,\nu,X))$ into itself.
This continuity implies that for any $t\in[0,1]$ the random variables
$\fhs _n \circ \Phi_t\circ \xi_1$  and $\fhs _n \circ \Phi_t\circ \xi_2 $ have the same laws on
$L ^  p([0,1] ;L^p(Z,\nu ;X))$.
Moreover,  by \eqref{cont-haar-shift}-(ii)
$\fhs _n \circ \Phi_t\circ \xi_i  \to \Phi_t\circ \xi_i$ in
$\CM^p  ([0,1] \times\Omega_i  ;L^p(Z,\nu,X))$.
\noindent
Taking into account the assumption that 
$\CL aw (\xi_1,\eta_1)=\CL aw (\xi_2,\eta_2)$ on $L ^ p([0,1] ;L ^ p(Z ,\nu;E) )
\times \CMME $, we conclude
by Corollary A.8 of \cite{Brz+Haus_2009-unique}
that for any $(t_1,\ldots , t_m) \in {[0,1]}^{m}$
\begin{eqnarray}
 & &  \mathcal{L}aw \bigl( {I}_{{\Phi }_{{t}_{1}}\circ {\xi }_{1}, {\tilde{\eta }}_{1}} ({t}_{1}),
 {I}_{{\Phi }_{{t}_{2}}\circ {\xi }_{1}, {\tilde{\eta }}_{1}} ({t}_{2}), \ldots ,
{I}_{{\Phi }_{{t}_{m}}\circ {\xi }_{1}, {\tilde{\eta }}_{1}} ({t}_{m}), {\xi }_{1}, {\tilde{\eta }}_{1}  \bigr) \label{eq:unique} \\
& & = \mathcal{L}aw \bigl( {I}_{{\Phi }_{{t}_{1}}\circ {\xi }_{2}, {\tilde{\eta }}_{2}} ({t}_{1}),
 {I}_{{\Phi }_{{t}_{2}}\circ {\xi }_{2}, {\tilde{\eta }}_{2}} ({t}_{2}), \ldots ,
{I}_{{\Phi }_{{t}_{m}}\circ {\xi }_{2}, {\tilde{\eta }}_{2}} ({t}_{m}), {\xi }_{2}, {\tilde{\eta }}_{2} \bigr)  \nonumber
\end{eqnarray}
on ${X}^{m} \times L^p([0,1] ;L^p(Z,\nu;E) ) \times {M}_I([0,1]\times Z)$.

We have to show that the triplets  $(
u_1,\xi_1,\eta_1)$ and $(u_2,\xi_2,\eta_2)$
have the same laws on $\DD([0,1];X)\times L ^ p([0,1] ;L^p([0,1],\nu ;E) )\times \CMME$.
To proceed further, let us recall the definition of the
dyadic projection $\pi_n :\DD([0,1];X)\to \DD([0,1];X)$ of order $n\in\mathbb{N}$
 introduced  in  \eqref{dyadic}, i.e.
\begin{equation}\label{dyadic-nnn}
(\pi _n   x)(t) := \sum_{j=0}^{2^n-1} \,1_{(j 2 ^ {-n},(j+1) 2 ^ {-n}]}(t) \; x( j 2 ^ {-n}),\; t\in[0,1].
\end{equation}
Note, that for  $i=1,2$ and $u_{ \xi_i,\tilde\eta_i}$ defined in \eqref{def_u_i}  the following identity holds
$$
\lk[ \pi_n \circ \, u_{ \xi_i,\tilde\eta_i} \rk]
(t) =  \sum_{j=0}^{2^n-1} \,1_{(j2 ^ {-n},(j+1) 2 ^ {-n}]}(t) \; I_{\Phi_{ j2 ^ {-n}}\circ \xi_i,\tilde \eta_i},\quad  t\in[0,1].
$$
By \eqref{eq:unique}, it follows that
$\pi_n \circ\, u_{ \xi_1,{\tilde{\eta}}_1}$ and $\pi_n \circ\, {u} _{\xi_2,{\tilde{\eta}}_2}$ have the same laws  on $\DD([0,1];X)$.
Finally, since by  Proposition \ref{chap:jacod}
the family of dyadic projections $\{\pi_n, n\in\mathbb{N}\}$ converges pointwise to the identity,
i.e.\ ${d}_{0}(\pi_n x,x) \to 0 $ \footnote{The Prohorov metric ${d}_{0}$ is defined in Appendix
\ref{app:Skorohod-space} in formula \eqref{eqn-dlog}.},
 for all $x\in \DD([0,1];X)$,
we infer by Proposition 2.6 of \cite{Brz+Haus_2009-unique} that the laws of the triplets $(u_i,\xi_i,\eta_i)$, $i=1,2$, are equal on
\[\DD([0,1];X)\times L^p([0,1];L^p(Z,\nu ;E))\times \CMME .\]
Hence the proof of  Theorem \ref{thm-main} is  complete.
\end{proof}

\bigskip
Note that existence of c\`{a}dl\`{a}g modification of the stochastic convolution process is fundamental both in the formulation and in the proof of Theorem \ref{thm-main}.
We will now give some examples of the spaces $X$ and $B$ whose existence is assumed in Theorem
\ref{thm-main}.

\bigskip
For simplicity we may assume that $A^{-1}$ exists and is a bounded operator on $E$.
For every $\beta>0$, we consider the extrapolation space $E_{-\beta}$ defined as the completion of the space $E$ with respect to the norm $\vert A^{-\beta} \cdot\vert_{E}$. Note that $A^\beta$ extends to an isometry, still denoted by $A^\beta$,  between $E$ and $E_{-\beta}$.

\bigskip
\begin{lemma} \label{lem-E_-1_cadlag}
Let  ${(S(t)}_{t\ge 0 }$ be a ${\mathcal{C}}_{0}$-semigroup generated by $-A$.  Then the process $u$ has an ${E}_{-1}$-valued \cadlag modification.
\end{lemma}

\begin{proof}
The assertion follows from \cite[Lemma 3.3]{Brz+Haus+Zhu_2012} according to which the process $u$ satisfies the following equality.
\begin{equation}  \label{eqn-strong}
 A^{-1}u(t) = \int_0^t u(s)\, ds+ \int_0 ^ t \int_Z A^{-1} \xi(s,z)\, \tilde \eta(ds,dz), \quad t\in [0,T].
\end{equation}
Hence, since by \cite{Brz+Haus_2009} the second term above has an $E$-valued \cadlag modification, we infer that the process $u$ has an $E_{-1}$-valued \cadlag modification.
\end{proof}

Using Theorem \ref{thm-main}, Lemma \ref{lem-E_-1_cadlag} and Remark \ref{R:Stoch_conv_L^p}, we obtain the following Corollary.

\begin{corollary}
Suppose that assumptions of Theorem \ref{thm-main} are satisfied.
Then  $\CL aw( (u _1,\xi  _1,\eta_1))=\CL aw((u _2,\xi _2,\eta_2)) $ on
\[ (\DD([0,T] ,{E}_{-1})\cap L^p([0,T] ;E))\times  L^p([0,T];L^p(Z,\nu ;E))\times \CMM .\]
\end{corollary}

\begin{proof}
By Lemma \ref{lem-E_-1_cadlag} the processes ${u}_{i}$, $i=1,2$, have ${E}_{-1}$-valued modifications.
By Remark \ref{R:Stoch_conv_L^p} ${\p }_{i}$-almost all paths of ${u}_{i}$, $i=1,2$, are in ${L}^{p}([0,T];E)$. Thus the assertion follows directly from Theorem \ref{thm-main}.
\end{proof}


\subsection{Stochasic convolution with a contraction type semigroup} \label{S:Contr_semigroup}
In paper \cite{Brz+Haus+Zhu_2012} it  is proven that
in the case when the space $E$ satisfies some stronger assumption and the semigroup is of contraction type
then the stochastic convolution  process has an $E$-valued c\`{a}dl\`{a}g modification. The problem of the existence of a c\`{a}dl\`{a}g modification is closely related to appropriate  maximal inequalities. We recall now  some of these results.

In addition we assume that the space $E$ satisfies the following condition
\begin{assumption} (see \cite{Brz+Haus+Zhu_2012}) \label{assump-E}
 There exists an equivalent norm ${\Vert \cdot \Vert}_{E} $ on $E$ and $q \in [p,\infty )$ such that the function
$\phi : E \ni x \mapsto {\Vert x  \Vert }_{E}^{q} \in \mathbb{R}$ is of class ${\mathcal{C}}^{2}$ and there exist constants ${k}_{1},{k}_{2}$ such that for every $x \in E$,
$ \vert \phi^{\prime}(x)\vert \le {k}_{1} {\Vert x \Vert }_{E}^{q-1}$  and $\vert \phi^{\prime\prime}(x)\vert \le {k}_{2} {\Vert x \Vert }_{E}^{q-2}$.
\end{assumption}

\begin{theorem} \label{thm-cadlag_mod} (\cite[Corollary 4.3]{Brz+Haus+Zhu_2012})
If  $E$ is a separable Banach space of martingale type $p$
satisfying Assumption \ref{assump-E}  and a ${\mathcal{C}}_{0}$-semigroup $(S(t))_{t \geq 0}$ is of contraction type, then there exists an $E$-valued c\`{a}dl\`{a}g modification $\tilde{u}$ of
$u$ such that for some constant $C>0$ independent of $u$ and all $t \in [0,T]$
and $0< r \le p $,
$$
   \e \sup_{0\le s \le t } {\Vert\tilde{u}\Vert}_{E}^{r }
 \le C \e {\Bigl( \int_{0}^{t} \int_{Z}{\Vert\xi (s,z)\Vert}_{E}^{p} \nu (dz)ds \Bigr) }^{\frac{r }{p}}.
$$
\end{theorem}

\bigskip
In the sequel,  in the case the assertion of Theorem \ref{thm-cadlag_mod} holds true,  by the stochastic convolution process we mean its $E$-valued c\`{a}dl\`{a}g modification.

Now, we are ready to  state the next corollary to our  main result Theorem \ref{thm-main}.
\begin{corollary} \label{cor-main} In addition to the assumptions of Theorem \ref{thm-main} we assume that
 the Banach space $(E,\vert \cdot \vert_E)$  satisfies  Assumption  \ref{assump-E} and the semigroup
$\big(S(t)\big)_{t\geq 0}$ is of contraction type.
Then
$\CL aw( (u _1,\xi  _1,\eta_1))=\CL aw((u _2,\xi _2,\eta_2)) $ on
$$ (\DD([0,T] ,E)\cap L^p([0,T] ;E))\times  L^p([0,T];L^p(Z,\nu ;E))\times \CMM .$$
\end{corollary}

\begin{remark}\label{rem-generalizations}
The only reason we have assumed in Corollary \ref{cor-main} that $\big(S(t)\big)_{t\geq 0}$ is a contraction type semigroup is that we need Theorem \ref{thm-cadlag_mod} about the existence of
 an $E$-valued  \cadlag modification of a stochastic convolution processes.
Corollary \ref{cor-main}
 remains valid for such class of semigroups $\big(S(t)\big)_{t\geq 0}$ for which the conclusion of
Theorem \ref{thm-cadlag_mod} holds true.
\end{remark}

\section{Applications}\label{sec-appl}

\subsection{Equations with a drift}\label{subsec-appl-drift}

Similarly to the paper \cite{Brz+Haus_2009-unique}
let us consider the following problem
\begin{equation} \label{E:SPDE}
\quad\quad\left\{\begin{array}{ll}
   d u(t) + Au(t) \; dt
 = b(t)\, dt +\int_{Z}\xi(t;z)\, \tilde \eta(dt;dz),\quad t \in [0,T],
\\  u(0) = 0.
\end{array}
\right.
\end{equation}
Let us suppose that $A$ is a linear operator satisfying Assumption \ref{assumption-semigroup}.
Moreover, assume that
\begin{trivlist}
\item[\bf (A.1)\rm ] $\xi : [0,T] \times \Omega  \to {L}^{p}(Z,\nu ;E)$ is a progressively measurable process such that $\e \bigl[ \int_{0}^{T} \int_{Z}\vert \xi (s,z) \vert_E ^{p} \, \nu (dz) ds  \bigr] < \infty $,
\item[\bf (A.2)\rm ] $b:[0,T]\times \Omega \to E $ is a progressively measurable process such that $\p $-a.s. $b \in {L}^{p}([0,T];E)$.
\end{trivlist}

\begin{definition}
Let $\xi $ and $b$ be two processes satisfying assumptions (A.1) and (A.2), respectively.
 A \sl solution \it to problem \eqref{E:SPDE} is
an $E$-valued predictable process $u$ such that $\p $-a.s.,  $u(0)=0$  and for every $t \in [0,T]$
the following equality
$$
   u(t) = \int_{0}^{t} S(t-s) b(s) \, ds + \int_{0}^{t} \int_{Z} S(t-s) \xi (s,z) \, \tilde{\eta }(ds,dz)
$$
holds $\p $-a.s.
\end{definition}

\begin{lemma} \label{L:SPDE_exist}
Let us assume that $(E,\vert \cdot \vert_E)$ be a separable  Banach space of martingale type $p$, where $1<p\le 2$. Let us also assume that  Assumption \ref{assumption-semigroup} is satisfied.
Let $\xi $ and $b$ be two processes satisfying assumptions (A.1) and (A.2), respectively.
Then there exists a solution to problem \eqref{E:SPDE}
 with  $\p $-almost all paths in ${L}^{p}([0,T];E)$.
Moreover, the solution has an ${E}_{-1}$-valued c\`{a}dl\`{a}g modification.

If in addition the space $E$ satisfies Assumption  \ref{assump-E} and the semigroup ${(S(t))}_{t \ge 0 }$ is of contraction type  then this solution has an $E$-valued
c\`{a}dl\`{a}g modification.
\end{lemma}

\begin{proof}[Proof of Lemma \ref{L:SPDE_exist}]
Since $S(t)$ is ${\mathcal{C}}_{0}$-semigroup and $b$ satisfies assumption (A.2), the process
$$
    v(t) = \int_{0}^{t} S(t-s) b(s) \, ds , \quad t \in [0,T]
$$
is well defined and $\p $-a.s. $v \in \mathcal{C}([0,T];E)$.
Thus the existence of a solution follows from Remark \ref{R:Stoch_conv_L^p} and   \cite[Lemma 3.1]{Brz+Haus+Zhu_2012}.
By Lemma \ref{lem-E_-1_cadlag} the solution has an ${E}_{-1}$-valued \cadlag modification.

The existence of an $E$-valued  c\`{a}dl\`{a}g modification
in the case when Assumption  \ref{assump-E} is satisfied and the semigroup ${(S(t))}_{t \ge 0 }$ is of contraction type
 follows directly from Theorem
\ref{thm-cadlag_mod}. The proof of the Lemma is thus complete.
\end{proof}

By methods  similar  to those used in the proofs of Theorem \ref{thm-main} and Lemma \ref{L:SPDE_exist} the following Corollary can be established.

\begin{corollary}\label{bbbb-sol}
Suppose that  the assumptions of Theorem \ref{thm-main} are satisfied.
Let $T>0$ and
 let  $u_i$, $i=1,2$, be the solutions of the following  problems
\begin{equation} \label{eqn-Langevin-Poisson-nnn}
\left\{\begin{array}{l}
d u_i(t) + Au_i(t  ) \; dt
= b_i(t)\; dt +\int_{Z}\xi_i(t, x)\, \tilde \eta_i(dt;dz),  \quad t\in [0,T],
\\  u_i(0) = 0  ,
\end{array}
\right.
\end{equation}
where $b_i:[0,T] \times\Omega_i  \to E $,
 $i=1,2$, is a progressively measurable process such that $ \PP$-a.s.\ $b_i\in L ^ p([0,T];E)$.
Assume that $\CL aw( (b_1,\xi  _1,\eta_1))=\CL aw((b_2,\xi _2,\eta_2)) $ on
$$
  L^p([0,T],E)\times L^p([0,T];L^p(Z,\nu;E))\times \CMM .
$$
Then
$\CL aw( (u _1,{b}_{1},\xi  _1,\eta_1))=\CL aw((u _2,{b}_{2},\xi _2,\eta_2)) $ on
\begin{eqnarray*}
& & (\DD([0,T] ,{E}_{-1})\cap L^p([0,T] ;E)) \\
& & \times {L}^{p}([0,T];E) \times  L^p([0,T];L^p(Z,\nu ;E))\times \CMM .
\end{eqnarray*}
\end{corollary}

\begin{proof}[Proof of Corollary \ref{bbbb-sol}]
The assertion follows directly from Lemma \ref{L:SPDE_exist} and Theorem \ref{thm-main}.
\end{proof}

\subsection{Equations with a contraction type semigroup}\label{subsec-appl-contraction}

\begin{corollary}\label{cor_eq_contr_semigroup}
If in addition to hypotheses of Corrolary \ref{bbbb-sol} the space $E$ satisfies Assumption \ref{assump-E}
and the semigroup ${(S(t))}_{t \ge 0}$ is of contraction type then
\begin{trivlist}
\item[(i)] the solution ${u}_{1}$ and ${u}_{2}$  have  $E$-valued
c\`{a}dl\`{a}g modifications,
\item[(ii)]
$\CL aw( (u _1,b_1,\xi  _1,\eta_1))=\CL aw((u _2,b_2,\xi _2,\eta_2)) $ on
$$\mathbb{D}([0,T],E) \times {L}^{p}([0,T],E)
\times  L^p([0,T];L^p(Z,\nu;E))\times \CMM.
$$
\end{trivlist}
\end{corollary}

\begin{proof}[Proof of Corollary \ref{cor_eq_contr_semigroup}]
The assertion is a direct consequence of Lemma \ref{L:SPDE_exist}
and  Theorem \ref{thm-main}.
\end{proof}

\subsection{Equations with an analytic semigroup}\label{subsec-appl-analytic}
In the last subsection we will make a stronger assumption than before.
\begin{assumption}\label{assumption-semigroup-analytic}
The operator $-A$ is an infinitesimal generator of  an analytic  ${\mathcal{C} }_{0}$-semigroup ${(S(t))}_{0\le t < \infty }$  on the space $E$.
\end{assumption}

We begin with the following useful result. An alternative proof of it could be obtained by applying \cite{Brz+Haus_2009} but for the convenience of the reader we have decided to include a self-contained proof.

\begin{lemma}\label{lem-estimates-analytic} 
 Suppose that that $(E,\vert \cdot \vert_E)$ is a separable  Banach space of martingale type $p$, where $1<p\le 2$. Suppose   also  that  Assumption \ref{assumption-semigroup} is satisfied. Then for every
 $ \alpha  \in (0, \frac{1}{p})$ there exists a constant ${C}_{\alpha ,p}>0$ such that  the stochastic convolution process $u$ defined  by
\[
    u(t):= \int_{0}^{t} \int_{Z} S(t-s) \xi (s,z) \tilde{\eta }(ds,dz), \qquad t \in [0,T].
\]
satisfies the following inequality
\begin{eqnarray*}
  \e \bigl[ {\| {A}^{\alpha }u \| }_{{L}^{p}(0,T;E)}^{p}\bigr]
 \le {C}_{\alpha ,p} \e \Bigl[ \int_{0}^{T} \int_{Z}\vert \xi (s,z) \vert_E^{p} \, d\nu (z)  ds \Bigr].
\end{eqnarray*}
\end{lemma}

\begin{proof}
Let us fix $ \alpha  \in (0, \frac{1}{p})$. Then, see \cite{Pazy_1983},
there exists a constant $C={C}(\alpha,T)>0$ such that for every $t>0$  linear operator ${A}^{\alpha }S(t)$ is well defined and bounded,  and
\[
    {|{A}^{\alpha }S(t)|}_{\mathcal{L}(E)} \le \frac{C }{{t}^{\alpha }},  \qquad t \in (0,T] .
\]

\bigskip  \noindent
We will show that for almost all $t\in [0,T] $
\begin{equation} \label{ineq-aux}
   \e \Bigl[  \int_{0}^{t}\int_{Z} \frac{1}{{(t-s)}^{\alpha p}} \vert \xi (s,z)  \vert_E^{p} \, \, \nu (dz)\, ds  \Bigr] < \infty .
\end{equation}
To this end it is sufficient to prove that
\[
 \int_{0}^{T} \Bigl\{  \e \Bigl[  \int_{0}^{t}\int_{Z} \frac{1}{{(t-s)}^{\alpha p}} \vert \xi (s,z)  \vert_E^{p} \, \, \nu (dz)\, ds  \Bigr]  \Bigr\} \, dt < \infty .
\]
Using the Fubini Theorem (for non-negative functions)  we obtain
\begin{eqnarray*}
& & \int_{0}^{T} \Bigl\{  \e \Bigl[  \int_{0}^{t}\int_{Z} \frac{1}{{(t-s)}^{\alpha p}} \vert \xi (s,z)  \vert_E^{p} \, \, \nu (dz)\, ds  \Bigr]  \Bigr\} \, dt \\
 & &=  \e \Bigl[  \int_{0}^{T} \Bigl\{  \int_{0}^{t}\frac{1}{{(t-s)}^{\alpha p}} \int_{Z}  \vert \xi (s,z)  \vert_E^{p} \, \, \nu (dz)\, ds  \Bigr\} \, dt  \Bigr] \\
& &  \le \e \Bigl[ \int_{0}^{T} \frac{1}{{s}^{\alpha p}} \, ds
  \cdot \int_{0}^{T} \int_{Z}  \vert \xi (s,z)  \vert_E^{p} \, \, \nu (dz)\, ds   \Bigr] \\
& &  = \frac{{T}^{-\alpha p +1}}{-\alpha p +1}
  \cdot  \e \Bigl[ \int_{0}^{T} \int_{Z}  \vert \xi (s,z)  \vert_E^{p} \, \, \nu (dz)\, ds   \Bigr] <\infty .
\end{eqnarray*}
Thus \eqref{ineq-aux} holds.
By \eqref{ineq-aux} we have for almost all  $t \in [0,T]$
\begin{eqnarray*}
  \e \bigl[ \vert  {A}^{\alpha }u(t) \vert_E ^{p}\bigr]
   & = & \e \Bigl[ \Bigl\vert \int_{0}^{t} \int_{Z} {A}^{\alpha }S(t-s)\xi (s,z)\tilde{\eta }(ds,dz)
  \Bigr\vert_E ^{p} \Bigr]  \\
   & \le  & \e \Bigl[  \int_{0}^{t}\int_{Z} \vert  {A}^{\alpha }S(t-s) \xi (s,z) \vert_E^{p} \, \, \nu (dz)\, ds  \Bigr] \\
 & \le  & C^{p}\e \Bigl[  \int_{0}^{t}\int_{Z} \frac{1}{{(t-s)}^{\alpha p}} \vert \xi (s,z)  \vert_E^{p} \, \, \nu (dz)\, ds  \Bigr] .
\end{eqnarray*}
Hence by the Fubini Theorem  we obtain
\begin{eqnarray*}
 \e \bigl[ \vert  {A}^{\alpha }u \vert_{{L}^{p}(0,T;E)}^{p} \bigr]
 & = & \e \Bigl[ \int_{0}^{T} \vert  {A}^{\alpha }u(t) \vert_E ^{p} \, dt \Bigr]
  =  \int_{0}^{T} \e \bigl[  \vert {A}^{\alpha }u(t) \vert ^{p} \,  \bigr] \, dt \\
 & \le & C^{p} \, \int_{0}^{T} \e \Bigl[ \int_{0}^{t}\int_{Z} \frac{1}{{(t-s)}^{\alpha p}} \vert \xi (s,z)  \vert_E^{p} \, \, \nu (dz)\, ds \Bigr] \, dt \\
 &\le & {C}_{\alpha ,p} \e \Bigl[ \int_{0}^{T} \int_{Z}\vert \xi (s,z) \vert_E^{p} \, d\nu (z)  ds \Bigr] ,
\end{eqnarray*}
where ${C}_{\alpha ,p}>0$ is some constant.

\end{proof}

We continue with the following strengthenings  of Lemmas \ref{lem-E_-1_cadlag} and \ref{L:SPDE_exist}.

\bigskip
\begin{lemma} \label{lem-analytic-cadlag} Under the assumptions of Lemma \ref{lem-estimates-analytic}
 the process $u$ has an ${E}_{\alpha-1}$-valued modification.
\end{lemma}

\begin{proof}From identity \ref{eqn-strong} we get
\begin{equation}  \label{eqn-strong-analytic}
 A^{\alpha-1}u(t) = \int_0^t A^\alpha u(s)\, ds+ \int_0 ^ t \int_Z A^{\alpha-1} \xi(s,z)\, \tilde \eta(ds,dz), \quad t\in [0,T].
\end{equation}
Hence, since by \cite{Brz+Haus_2009} the second term above has an $E$-valued \cadlag modification and by the previous Lemma \ref{lem-estimates-analytic} the first term above has an $E$-valued continuous modification,   we infer that the process $u$ has an $E_{\alpha-1}$-valued \cadlag modification.
\end{proof}

\begin{lemma} \label{lem-SPDE_exist-analytic}
Assume that $E$ be a separable  Banach space of martingale type $p$, $1<p\le 2$,
that Assumption \ref{assumption-semigroup-analytic} holds and that $\alpha < \frac1p$.
Let $\xi $ and $b$ be two processes satisfying assumptions (A.1) and (A.2), respectively.
Then there exist a solution $u$ to problem \eqref{E:SPDE} with  $\p $-almost all paths in ${L}^{p}([0,T];D(A^\alpha))$.
Moreover, every  this solution has an ${E}_{\alpha-1}$-valued c\`{a}dl\`{a}g modification.
\end{lemma}
\begin{proof}The proof is obvious.
\end{proof}
We finish this subsection with the main result of it.
\begin{corollary}\label{cor-main-analytic}
Suppose that  the assumptions of Theorem \ref{thm-main} are satisfied and that Assumption \ref{assumption-semigroup-analytic} holds. Assume also that $\alpha\in (0,\frac1p)$.
\noindent
Let $T>0$ and
 let  $u_i$, $i=1,2$, be the solutions of the following  problems
\begin{equation} \label{eqn-Langevin-Poisson-nnn-analytic}
\left\{\begin{array}{l}
d u_i(t) + Au_i(t  ) \; dt
= b_i(t)\; dt +\int_{Z}\xi_i(t, x)\, \tilde \eta_i(dt;dz),  \quad t\in [0,T],
\\  u_i(0) = 0  ,
\end{array}
\right.
\end{equation}
where $b_i:[0,T] \times\Omega_i  \to E $,
 $i=1,2$, is a progressively measurable process such that $ \PP$-a.s.\ $b_i\in L ^ p([0,T];E)$.
Assume that $\CL aw( (b_1,\xi  _1,\eta_1))=\CL aw((b_2,\xi _2,\eta_2)) $ on
$$
  L^p([0,T],E)\times L^p([0,T];L^p(Z,\nu;E))\times \CMM .
$$
Then
$\CL aw( (u _1,{b}_{1},\xi  _1,\eta_1))=\CL aw((u _2,{b}_{2},\xi _2,\eta_2)) $ on
\begin{eqnarray*}
& & (\DD([0,T] ,{E}_{\alpha-1})\cap L^p([0,T] ;D(A^\alpha))) \\
& & \times {L}^{p}([0,T];E) \times  L^p([0,T];L^p(Z,\nu ;E))\times \CMM .
\end{eqnarray*}

\end{corollary}

\appendix

\section{Equality of Laws}\label{app:A}

The purpose of this Appendix is to prove the following essential result. This result should be known but we were  unable to trace it in the literature.
\begin{proposition}\label{prop-main} 
Assume that $X$  and $E$ are Polish spaces such that $E\subset X$ and $E\subset X$  and the natural embedding
$\iota :E\hookrightarrow  X$ is  continuous.
For fixed $k=1,2$, let
${\mathfrak{A}}_{k} =(\Omega_k,\CF_k,\PP_k)$ be a probability space, and
$\xi _k$ be an  $E$-valued random variable
defined on ${\mathfrak{A}}_{k}$.
If $\CLaw(\iota\circ \xi_1)=\CLaw(\iota\circ \xi_2)$ on $X$, then $\CLaw(\xi_1)=\CLaw(\xi_2)$ on $E$.
\end{proposition}

\begin{proof}
It is enough to prove that $\mathcal{B}(E)=\iota^{-1}(\mathcal{B}(X))$.
Since $\iota$ is a continuous mapping,
the inclusion $\supset$ follows. To prove the inclusion $\subset$ let us take $A\in\mathcal{B}(E)$.  Therefore, by the Kuratowski Theorem \cite[p. 499]{Kuratowski_1952}, see also
\cite[Theorem 1.1, p. 5]{Vakh+Tar_1987} and \cite[\S 1.3]{para}, since $\iota$ is an injective continuous (and hence Borel)
mapping, the set $\iota(A)$ belongs to $\mathcal{B}(X)$. By injectivity of $\iota$ we infer that
$A=\iota^{-1}\big(\iota(A)\big)$ which completes the proof.
\end{proof}

\del{We can also give a more direct proof of Proposition
\ref{prop-main} that does not recall to the deep result of
Kuratowski.

\begin{proof}
Let $A\in\CB(E)$. Let us fix  $\ep>0$. Since  $E$ is Polish by   \cite[Theorem II.3.2]{para}
it follows that the laws of both processes $\xi_k$, $k=1,2$,  are
Radon measures on $(E, \mathcal{B}(E))$. Hence
 there exists a compact subset $K$ of  $E$ 
such that $K\subset A$ and  $\PP_1( \xi_1\in A\cap K^c  
) +\PP_2( \xi_2\in A\cap K^c 
)<\ep$.
Next, by  the triangle inequality  we have
\DEQS
\lqq{
|\PP_1(\xi_1\in A)-\PP_2(\xi_2\in A)|
 }
\\
&\le & \lk|\PP_1(\xi_1\in A\cap K^c )+ \PP_1(\xi_1\in K)- \PP_2(\xi_2\in A\cap K^c 
 )- \PP_2(\xi_2\in K)\rk|
\\
&\le &\PP_1(\xi_1\in A\cap K^c 
 )+  \PP_2(\xi_2\in A\cap K^c 
 ) + \lk| \PP_1(\xi_1\in K) - \PP_2(\xi_2\in K)\rk|
\\
&\le &\ep +  \lk| \PP_1(\xi_1\in K) - \PP_2(\xi_2\in K)\rk|.
\EEQS
Since embedding $\iota:E\hookrightarrow X$ is continuous, the set $\iota(K)$ is  compact in $X$ and hence
a Borel subset of $X$.
Since $\PP_k(\xi_k\in K)=\PP_k(\iota \circ \xi_k\in \iota (K))$, for $k=1,2$ and by assumptions
$\PP_1(\iota \circ \xi_1\in \iota (K))= \PP_2(\iota \circ \xi_2\in \iota (K))$ we infer that $\PP_1(\xi_1\in K)=
\PP_2(\xi_2\in K)$. This implies that $| \PP_1(\xi_1\in A)-\PP_2(\xi_2\in A)|< \ep$, which completes the proof.
\end{proof}
}

\begin{corollary}\label{cor-prop-main}
Assume that $X$, $E$ and $F$ are Polish spaces such that $E\subset X$ and $F\subset X$,  and the natural embeddings
$\iota_E :E\hookrightarrow  X$ and $\iota_F :F\hookrightarrow  X$ are continuous. If the random variables $\xi_i$, $i=1,2$ as in Proposition \ref{prop-main} are $E \cap F$ valued and $\CLaw(\xi_1)=\CLaw( \xi_2)$ on $E$ then  $\CLaw(\xi_1)=\CLaw( \xi_2)$ on $F$.
\end{corollary}
\begin{proof}[Proof of Corollary \ref{cor-prop-main}] By assumptions, $\CLaw(\iota_E\circ \xi_1)=\CLaw(\iota_E\circ \xi_2)$ on $X$. But
$\iota_E\circ \xi_i=\iota_F\circ \xi_i$, $i=1,2$ and hence $\CLaw(\iota_F\circ \xi_1)=\CLaw(\iota_F\circ \xi_2)$ on $X$. Applying next Proposition
\ref{prop-main} concludes the proof.
\end{proof}

\section{$L^p$-spaces and the Haar system}\label{app:haar-system}

In this Appendix we assume that $(Y,|\cdot |)$ is a separable Banach space.
As in the Proof  of Theorem \ref{thm-main} we assume for simplicity that $T=1$.
In this section we recall some facts about the approximation properties of the Haar system.
For $n\in\mathbb{N}$, let  $\Pi^n=\{ s^n_0=0<s^n_1<\cdots <s^n_{2^n}\}$
be a  partition of the interval $[0,1]$ defined  by $s_j ^ n=j\,2 ^ {-n}$, $j=1,\cdots, 2^n$.
Each interval of the form $(s_{j-1} ^ n,s_j ^ n]$, where $n\in \mathbb{N}$ and $j=1,\ldots, 2^n$ is called a \it dyadic \rm interval.
For  $n\in \mathbb{N}$, the \it $j^ {th}$ element, for $j=1,\ldots, 2^n$,  of the Haar system of order $n$ \rm is the indicator function of the interval $(s_{j-1}^n,s_j^{n}]$, i.e.
$1_{(s_{j-1}{n},s_j^{n}]}$.
For $n\in\mathbb{N}$ let $\fh_n: L ^ p([0,1],Y)\to L ^ p([0,1],Y)$ be the \it Haar projection \rm of order $n$, i.e.
\begin{equation} \label{haar-projection}
\fh_n ( x) = \sum_{j=1} ^{2^n} 1_{(s^n_{j-1},s^n_{j}]} \otimes \iota_{j,n}(x), \;\; x\in  L ^ p([0,1],Y),
\end{equation}
where $\iota_{j,n}:  L ^ p([0,1],Y)\to Y$,
is the averaging operator over the interval $(s^n_{j-1},s^n_{j}]$, i.e.
\begin{equation}  \label{xdef-ij}
 \iota _{j,n}(x) :=
\frac{1}{s^n_{j}-s^n_{j-1}}\; \int_{s^ n_{j-1}}^{s ^ n_{j}} x(s)\,ds, \;  x\in  L ^ p([0,1],Y).
\end{equation}
In the above, for $f\in L ^ p([0,1],\mathbb{R})$ and $y\in  Y$, by $f\otimes y$ we mean an element of $L ^ p([0,1],Y)$ defined by
$[0,1]\ni t\mapsto f(t) y\in Y$.

\begin{remark}\label{domination-haar} Note that by the Jensen inequality, every  map
$\iota_{j,n}:  L ^ p([0,1],Y)\to Y$ is a linear contraction. Therefore, since $ \vert 1_{(s^n_{j-1},s^n_{j}]}\vert_{L ^ p([0,1],\mathbb{R})}=
2^{-n}$, we infer that that for any $n\in\mathbb{N}$
 \begin{equation}
 \label{eqn-contraction}
  \|\fh _n (x) \|_{ L ^ p([0,1],Y)}\leq \| x\|_{L ^ p([0,1],Y)}, \quad   x\in L ^ p([0,1],Y).
\end{equation}

Therefore,  since $\fh_n x\to x$ in $L ^ p([0,1],Y)$ for any $x\in C^1([0,1],Y)$ and $C^1([0,1],Y)$ is a dense subspace of  $L ^ p([0,1];Y)$, we infer that
\DEQSZ \label{haar-conv}
 \fh _n (x) \to x \; \mbox{ in } \; L^p([0,1];Y),
\quad \forall x\in  L ^ p([0,1];Y).
\EEQSZ
\end{remark}

Let $\xi$ be a $Y$-valued progressively measurable $p$ integrable  stochastic process. One way to get a sequence of step functions is to approximate $\xi$ by the sequence $(\fh_n \circ \xi)_{n\in\mathbb{N}}$, where
$\fh_n$ is the Haar  projection of order  $n$. The only problem which arises is,
 that    $\fh_n \circ \xi$ is not necessarily progressively measurable. To get this property, we have to shift the projection by one time interval. With $\iota_{j,n}$
being defined in  \eqref{xdef-ij}, the $n$-th order shifted Haar projection is a linear bounded map  in  $ L ^ p([0,1];Y)$  is defined by

\DEQSZ\label{haar-projection-shifted}
\fhs_n  x =\sum_{j=1} ^{2^n -1} 1_{(s^n_j,s^n_{j+1}]}\otimes \iota_{j-1,n}(x), \;\; x\in L ^ p([0,1];Y),
\EEQSZ
where we put ${\iota }_{0,n}=0$ for every $n\in\mathbb{N}$.
\begin{proposition}
\label{cont-haar-shift}
The following holds:
\begin{trivlist}
\item[(i)]
For any $n\in\mathbb{N}$, the shifted Haar projection $\fhs_n:L ^ p([0,1];Y)\to L ^ p([0,1];Y)$ is a continuous operator.
\item[(ii)]
For all $x\in  L ^ p([0,1];Y)$, $\fhs_n x\to x$ in $L ^ p([0,1];Y)$.
\end{trivlist} 
\end{proposition}
\begin{proof}[Proof of Proposition \ref{cont-haar-shift}]
For each $n\in {\mathbb{N}}$ let us define
$$ (\fs_n x)(s) := \bcase x(0),& \mbox{ if } s\le 2 ^ {-n},\\
 x(s-2 ^ {-n}),&  \mbox{ if } s\in ( 2 ^ {-n},1],
\ecase
$$
Obviously, $\fs_n$ is a linear contraction in ${L}^{p}([0,1];Y)$. Moreover,
$$
      \fhs_n = \mathfrak{h} \circ \fs_n.
$$
Therefore part (i) follows from Remark \ref{domination-haar}.
Since $\fs_n x\to x$  for all $x\in C^1([0,1],Y)$, assertion (ii) follows by similar arguments as in
\eqref{haar-conv}. This completes the proof.
\end{proof}

\section{The Skorokhod space}\label{app:Skorohod-space}

For an introduction to the Skorokhod space we refer to Billingsley \cite{billingsley}, Ethier and
Kurtz \cite{MR838085} and Jacod and Shiryaev
\cite{jacod}. In this Section we state only these results which are necessary for our work.

Let $(Y,{|\cdot |}_{Y})$ be a separable Banach space.
The space $\DD([0,1];Y)$ denotes the space of all right continuous functions $x:[0,1]\to Y$
with left-hand limits. Let $\Lambda $ denote the class of all strictly increasing continuous functions $\lambda:[0,1]\to [0,1]$ such that $\lambda(0)=0$ and {$\lambda(1)=1$}.
Obviously any element $\lambda \in \Lambda$ is a homeomorphism of $[0,1]$ onto itself.
Let us define the Prohorov metric $d_0$ by
\begin{equation}\label{eqn-dlog}
\begin{array}{rcl}
\Vert \lambda\Vert_{\rm log} &:=& \sup_{t\not=s\in [0,1]} \Big\vert \log \frac{\lambda(t)-\lambda(s)}{t-s}\Big\vert,\\
\Lambda_{\rm log}&:=& \big\{ \lambda \in \Lambda:  \Vert \lambda\Vert_{\rm log} <\infty\big\}\\
d_0(x,y)&:=&\inf \Big\{ \Vert \lambda\Vert_{\rm log} \vee  \sup_{t\in [0,1]} {|x(t)-y(\lambda(t))|}_{Y}: \, \lambda \in \Lambda_{\rm log}\Big\}.
\end{array}
\end{equation}
The space $\DD([0,1];Y)$ equipped with the metric $d_0$ is a separable complete metric space.
We recall the notion of the dyadic projection.
\begin{definition} Assume that $n\in \mathbb{N}$.
The $n$-th order  dyadic projection is a linear map $\pi_n: \DD([0,1];Y)\to\DD([0,1];Y)$ defined by
\DEQSZ\label{dyadic}
& & \pi_n x :=
 \sum_{i=0}^{2 ^ n-1} \,1_{(2 ^ {-n}i,2 ^ {-n}(i+1)]} \, x( 2 ^ {-n}i), \quad x\in \DD([0,1];Y) .
\EEQSZ
\end{definition}
An  important  property  of  the dyadic projection is  given
in the following result.

\begin{proposition}\label{chap:jacod} (see \cite[Proposition B.5]{Brz+Haus_2009-unique})
If  $x\in \DD([0,1];Y)$ then
$$
    \lim_{n \to \infty } {d}_{0}(x,{\pi }_{n}x) =0 .
$$
\end{proposition}

\def\cprime{$'$} \def\cprime{$'$} \def\cprime{$'$} \def\cprime{$'$}

\end{document}